\documentclass{amsart}

\usepackage{amsmath}
\usepackage{amsfonts,amssymb,setspace}

\begin{document}

\newtheorem{theorem}{Theorem}
\newtheorem{corollary}[theorem]{Corollary}
\newtheorem{thm}{Theorem}[section]
\newtheorem{cjct}{Conjecture}
\newtheorem{cor}[thm]{Corollary}
\newtheorem{lem}[thm]{Lemma}
\newtheorem{prop}[thm]{Proposition}
\theoremstyle{definition}
\newtheorem{defn}[thm]{Definition}
\theoremstyle{remark}
\newtheorem{rem}[thm]{Remark}
\newtheorem{examp}{Example}
\numberwithin{equation}{section}

\newcommand{\R}{\mathbb{R}}
\newcommand{\ZM}{\mathbb{Z}}
\newcommand{\QM}{\mathbb{Q}}
\newcommand{\NM}{\mathbb{N}}
\newcommand{\CM}{\mathbb{C}}

\newcommand{\eps}{\varepsilon}

\newcommand{\Z}{\mathbb{Z}}
\newcommand{\conj}[1]{\overline#1}
\newcommand{\diag}{\rm diag}
\newcommand{\disps}{\displaystyle}
\newcommand{\ds}{\displaystyle}
\newcommand{\dsum}{\ds\sum}
\newcommand{\ep}{\qed\endtrivlist}
\newcommand{\eqskip}{ \vspace*{2mm} }
\newcommand{\dint}{ {\disps \int} }
\newcommand{\kk}{\mathcal{K}}
\newcommand{\ko}{K_{0}}
\newcommand{\mr}{ \mathcal{R} }
\newcommand{\mth}{ \Theta }
\newcommand{\nsu}{ \mathcal{N}_{u} }
\newcommand{\om}{ \omega }
\newcommand{\parag}{ \hspace{-0,65cm} }
\newcommand{\sca}{ \mathcal{S} }
\newcommand{\ub}{\bar{u}}
\newcommand{\vphi}{\varphi}


\title[Zero--free half--planes of Riemann's zeta 
function]
{A Li--type criterion for zero--free half--planes
of Riemann's zeta function }
\author{Pedro Freitas}



\address{Departamento de Matem\'{a}tica,
Instituto Superior T\'{e}cnico, Av.Rovisco Pais, 1049-001 Lisboa,
Portugal.} \email{pfreitas@math.ist.utl.pt}
\date{\today}

\begin{abstract}
    We define a sequence of real functions which coincide 
    with Li's coefficients at one and which allow us to
    extend Li's criterion for the Riemann Hypothesis
    to yield a necessary and sufficient condition for the existence 
    of zero--free strips inside the critical strip $0<\Re(z)<1$. 
    We study some of the properties of these functions, including their 
    oscillatory behaviour.
\end{abstract}

\maketitle
\section{Introduction}

In 1997 Li gave a necessary and sufficient condition for the Riemann 
Hypothesis (RH) to hold, based on the positivity of a sequence of real 
numbers~\cite{li}. Li's coefficients may be written as
\begin{equation}\label{liscoeff}
\lambda_{n+1}= \frac{\ds 1}{\ds n!} \frac{\ds d^{n+1}}{\ds d 
s^{n+1}}\left[ s^n \log\left(\xi(s)\right)\right]_{s=1}, \;\; n=0,1,\ldots,
\end{equation}
and the criterion then states that the RH holds if and only if 
$\lambda_{n}$ is greater than or equal to zero for all $n=1,2,\ldots$. 
Here $\xi$ is the function defined by
\[
\xi(s)=s(s-1)\pi^{-s/2}\Gamma(s/2)\zeta(s),
\]
where $\Gamma$ and $\zeta$ are, respectively, Euler's gamma and Riemann's zeta
functions. As is well known, $\xi$ is an entire function satisfying the 
functional equation $\xi(s)=\xi(1-s)$. 

Equivalently, these coefficients may also be defined as the coefficients in
the Taylor series
of the function
\begin{equation}\label{taylorserphi}
\psi(z) = \frac{\ds \vphi'(z)}{\ds \vphi(z)} = \dsum_{n=0}^{\infty}\lambda_{n+1} 
z^n
\end{equation}
where $\vphi$ is given by
\begin{equation}\label{phidefn}
\vphi(z) = \xi\left(\frac{\ds 1}{\ds 1-z}\right).
\end{equation}
In fact, and as remarked in~\cite{bola}, the idea behind Li's criterion is that
the RH holds if and  only if the coefficients $\lambda_{n}$ satisfy a growth
condition, namely, that
$\lim\sup |\lambda_{n+1}|^{1/n}\leq1$, and that this may in fact be 
replaced by the (one--sided) condition on their positivity.
From this it follows that the positivity of Li's coefficients is related
to the derivatives of the functions $\xi$ and 
$\vphi$.

At this point it should be mentioned that a set of 
coefficients which differ from Li's only by a (positive) 
multiplicative constant was introduced earlier by Keiper 
in~\cite{keip}. In particular, it is already stated in~\cite{keip} 
that the positivity of these coefficients is implied by the RH. I am
indebted to Andr\'{e} Voros for pointing this
out to me.

Li's coefficients have recently been the subject of much study, from both
theoretical and numerical points of 
view~--~see~\cite{brow,bola,coff,coff2,laga,masl1,masl2,voro},
for instance. As was pointed out in~\cite{bola}, this criterion is 
not specific for the function $\xi$, and will work for any function which
can be written as a product of the form
\begin{equation}\label{prodf}
f(z) = \prod_{\rho\in Z}\left(1-\frac{\ds z}{\ds \rho}\right),
\end{equation}
or also in the more general case of a multi--set, provided some 
conditions are imposed on the asymptotic behaviour of the zero set 
$Z$.
The main purpose of the present note is to show that Li's criterion may be 
extended in a simple and natural way to provide a necessary and sufficient 
condition for the existence of a
zero--free half--plane  $\Re(z)>\tau/2$, and to study the behaviour of the
functions that arise naturally while carrying out
this extension. We remark that an extension of Li's criterion to 
provide zero--free regions for Riemann's zeta function has been 
considered in~\cite{brow}. However, this uses a completely different 
approach based on an analysis of the implications of the non--negativity 
of the first $N$ of Li's coefficients.

Although most of the results 
presented here apply to a more general function $f$ as above, we will 
concentrate on the case where $f=\xi$. In this case, and for $\tau$ 
in $(1,2)$, the condition that there are no zeros of zeta with 
$\Re(z)>\tau/2$ is sometimes referred to as the
quasi--RH, and we remark that this weaker version of the RH also remains an 
open problem. More precisely, it is not known if there exists 
$\tau_{0}$ on $(1/2,1)$ such that there are no zeros in the strip 
$\Re(s)>\tau_{0}$.

In order to extend Li's criterion, let us define the coefficients
\begin{equation}
    \label{dcoeff}
    \alpha_{n+1}(\tau) = \left.\frac{\ds 1}{\ds n!} \frac{\ds d^{n+1}}{\ds d 
s^{n+1}}\left[ s^{n}
\log\left(\xi(s)\right)\right]\right|_{s=\tau}, \;\; 
n=0,1,\ldots,
\end{equation}
where $\tau$ is a positive number. We remark that 
of the original (implicit) definitions of the coefficients $\lambda_{n}$
in~\cite{li}, not all correspond to the expression above. 
More precisely, apart from equations~(\ref{liscoeff}) 
and~(\ref{taylorserphi}), the expression
\[
\sum_{\rho\in Z}\left[1-\left(1-\frac{\ds 1}{\ds \rho}\right)^{n}\right]
\]
is also used in that paper as corresponding to $\lambda_{n}$, while in fact this 
would correspond to $\lambda_{-n}$ -- see~\cite{bola} and also 
Section~\ref{proofmain} below. In fact, at $\tau=1$ it is indeed true 
that $\lambda_{n}=\lambda_{-n}$, and thus this distinction was not important in
Li's paper. However, since this will not be the case for other values 
of $\tau$, amd we use~(\ref{dcoeff}) as a 
starting point, there is a discrepancy between the notation in our paper and
that used in other papers, such as~\cite{bola}.

Clearly these functions are 
analytic for real values of $\tau$, and have singularities at the 
zeros of the function $\xi$.
The main idea here is that while for $\tau$ equal to one we recover
Li's coefficients, the positivity of the $\alpha_{n}$'s for other values of 
$\tau$ is equivalent to the non--existence of zeros in the half--plane 
$\Re(s)>\tau/2$. In particular, and since for $\tau$ smaller than 
one we obtain a region which contains the critical line, in this 
case some of the coefficients $\alpha_{n}$ must take on negative values. On the other hand, 
for $\tau$ greater than or equal to two we are outside the critical 
strip and hence the $\alpha_{n}$'s must all be non--negative.

More precisely, we have the following
\begin{theorem}\label{mainthm}
    Given $\tau$ on $[1/2,\infty)$, the half--plane $\Re(s)>\tau/2$ is
    a zero--free region for the Riemann zeta
    function if and only if $\alpha_{n}(\tau)$ is nonnegative for all 
    positive integer $n$.
\end{theorem}

Note that Theorem 1 in~\cite{bola} already allows us to 
obtain a condition which will ensure that the zeros of the 
function $\xi$ are confined to a half--plane -- this may be done simply 
by replacing $\rho$ by $\rho/\tau$ in the relevant expressions. The 
main contribution of the above result is to relate this to the positivity 
of the coefficients $\alpha_{n}$ defined by~(\ref{dcoeff}). The proof 
also shows that these functions are directly related to the Taylor 
coefficients of the function $\psi$ around the point $z=1-1/\tau$.

Since we will follow an approach similar to that used in Li's paper, in order
to prove Theorem~\ref{mainthm} we will
need some results concerning Taylor expansions of
the function $\vphi(z)=\xi(1/(1-z))$. This will be done in 
Section~\ref{cnprop}, where the behaviour of the corresponding Taylor
coefficients is analysed. We shall present some simpler
proofs for some known results, and also obtain some new monotonicity 
properties which we believe to be of interest in their own right.

In Section~\ref{proofmain} we prove Theorem~\ref{mainthm}, and 
in Section~\ref{prop} we obtain some further properties of 
the coefficients $\alpha_{n}$. These include a system of differential 
equations satisfied by the $\alpha_{n}$'s -- see 
Theorem~\ref{tsysdiff} -- and some oscillatory properties which are
summarized in the following
\begin{theorem}\label{oscilthm}
    For $n=1,2,\ldots$ and $0\leq\tau<2$ the following facts hold
    \begin{itemize}
	\item[(i)] $ \alpha_{n}(0) = n\frac{\ds \xi'(0)}{\ds \xi(0)}<0$;
	\item[(ii)] there exists a sequence $1/2=a_{1}>a_{2}>\ldots$ such that
	$\alpha_{n}(a_{n})=0$ and $\alpha_{n}(\tau)<0$ for $\tau<a_{n}$;
	\item[(iii)] between any two real zeros of $\alpha_{n}(\tau)$ there
	exists at least 
	one real zero of $\alpha_{n+1}(\tau)$;
	\item[(iv)] for any positive integer $N$ there exists a positive 
	integer $n_{0}$ such that the function $\alpha_{n}(\tau)$ has at least
	$N$ zeros on the interval $(0,2)$, for all $n$ larger than $n_{0}$.
    \end{itemize}
\end{theorem}

We then proceed to extend the definition of the coefficients 
$\alpha_{n}$ to the case of real $n$. This is done in 
Section~\ref{secf}, where some alternative representations for this
more general function are given. 
Finally, in the last section, we briefly discuss the results obtained.

\section{\label{cnprop}Some properties of the Taylor coefficients of 
the functions~$\vphi$~and~$\xi$ }

We begin by proving a simple lemma which will be useful 
in the sequel. This can easilly be obtained from known results 
(see~\cite{knop}, for instance), but we include it here in a form 
appropriate for our purposes.
\begin{lem}\label{simplelem}
    Let $f$ be an analytic function in an open subset $\Omega$ of 
    $\CM$, and let $z_{0}<z_{1}$ be two (real) points in $\Omega$. 
    Write
    \[
    f(z) = \dsum_{n=0}^{\infty}a_{n}(z-z_{0})^n = 
    \dsum_{n=0}^{\infty}b_{n}(z-z_{1})^n,
    \]
    and assume the following conditions hold
    \begin{enumerate}
	\item[(C1)] $a_{n}$ is real and non-negative for $n=0,1,\ldots$
	\item[(C2)] the first series is convergent at $z=z_{1}$
    \end{enumerate}
    Then $b_{n}\geq0$ for $n=0,1,\ldots$.
    Furthermore, if $a_{m}>0$ for some integer $m$, then $b_{n}>0$ 
    for $n=0,\ldots,m$.
\end{lem}
\begin{proof}
    From
    \[
    f(z) = \dsum_{n=0}^{\infty}a_{n}(z-z_{0})^n = 
    \dsum_{n=0}^{\infty}b_{n}(z-z_{1})^n,
    \]
    and taking $p$ derivatives with respect to $z$, we have that
    \[
    \dsum_{n=p}^{\infty}\left(\begin{array}{c}n\\ 
    p\end{array}\right)
    a_{n}(z-z_{0})^{n-p} = 
    \dsum_{n=p}^{\infty}\left(\begin{array}{c}n\\ 
    p\end{array}\right)b_{n}(z-z_{1})^{n-p}.
    \]
    Making now $z=z_{1}$ yields
    \[
    b_{p} = \dsum_{n=p}^{\infty}\left(\begin{array}{c}n\\ 
    p\end{array}\right)a_{n}
    (z_{1}-z_{0})^{n-p},
    \]
    from which the result follows.
\end{proof}

Using this we obtain a simple proof of the following result regarding 
the sign of the derivatives of the function $\xi$, which can 
already be found in~\cite{coff}. This is a generalization of a result 
of Pustyl'nikov's, who proved in~\cite{pust} that a necessary 
condition for the RH to hold is that
all even derivatives of $\xi$ at the point $z=1/2$ are positive. While 
these proofs are based on the integral representation of the 
function $\pi^{-z/2}\Gamma(z/2)\zeta(z)$, here we show this by
direct series manipulation and the result at $z=1/2$, which may
already be found on page 41 of~\cite{edwa}.
\begin{thm}\label{extension}
    Even derivatives of $\xi$ are positive for all real values of $z$, 
    while odd derivatives are positive for $z>1/2$ and negative for 
    $z<1/2$.
\end{thm}

\begin{proof} Consider the
series development of $\xi$ at a point $z_{0}$, that is,
\begin{equation}\label{xiseries}
\xi(z) = \dsum_{n=0}^{\infty} s_{n}(z_{0})(z-z_{0})^n.
\end{equation}
Since $\xi$ is an entire function,
\[
s_{0}(1/2)=\xi(1/2) = 
-\Gamma(1/4)\zeta(1/2)/(4\pi^{1/4})\approx .994242,
\]
and
we have that $s_{2n}(1/2)$ is strictly positive, it 
follows from the lemma that $s_{n}(z_{0})$ is strictly positive for 
all $n$ and all real $z$ larger than $1/2$. The result for 
$z$ less than $1/2$ now follows by repeated differentiation of the 
functional equation for $\xi$.
\end{proof}

Write now the series development for the function $\vphi$ around a real 
number $z_{0}$ as
\[
\vphi(z) = \dsum_{n=0}^{\infty} c_{n}(z_{0}) (z-z_{0})^n.
\]
This series is well defined for all $z_{0}$ different 
from one, and its radius of convergence is equal to
$|1-z_{0}|$.

We shall first show that the coefficients $c_{n}(-1)$ 
are all positive with the exception of $c_{1}(-1)$ which vanishes -- 
in what follows, and provided there is no room for confusion, we shall 
drop the explicit reference to the point $z_{0}$.

\begin{thm}\label{posit} For all positive integers $n$ different from 
$1$, we
    have that
    $\vphi^{(n)}(-1)>0$, while $\vphi'(-1)=0$.
\end{thm}
\begin{proof}
    From the definition of $\vphi$ we see that $\vphi(-1) = 
    \xi(1/2)>0$. Differentiating both sides of~(\ref{phidefn}) with respect
    to $z$ yields
    \[
    \vphi'(z) = 
    \xi'\left(\frac{\ds 1}{\ds 1-z}\right) \left(\frac{\ds 1}{\ds 
    1-z}\right)^{2}.
    \]
    At $z=-1$ this gives $\vphi'(-1)=0$, while further repeated differentiation
    of this identity will now yield that the 
    $k^{\rm th}$ derivative of $\vphi$ at $-1$ is a linear combination 
    with positive coefficients of
    derivatives of $\xi$ taken at $1/2$. Since for $k$ larger than 
    $1$ there will always exist at least one even derivative of 
    $\xi$, the result follows.
\end{proof}

An immediate consequence of this result is 
the positivity of the coefficients in the series development of 
$\vphi$ around any point on the interval $(-1,1)$.

\begin{cor}\label{cpos} 
    For all $z$ in $(-1,1)$, $c_{n}(z)>0$.
\end{cor}
\begin{proof}
    It follows by direct application of Lemma~\ref{simplelem} with 
    $z_{0}=-1$, taking 
    into account that $\vphi$'s only singularity is at $z=1$, and thus 
    condition (C2) in the lemma is satisfied, while (C1) is a direct 
    consequence of Theorem~\ref{posit}
\end{proof}
\begin{rem} The coefficients $c_{n}(0)$ correspond to the 
coefficients $a_{n}$ in Li's paper.
\end{rem}
\begin{rem} This may also be proven in the same way as in 
Theorem~\ref{posit}, using the positivity results of Theorem~\ref{extension}.
\end{rem}

As a consequence of the positivity of the coefficients $c_{n}$ we see that 
they are actually increasing functions of $z_{0}$.

\begin{cor} \label{monotz}
    The functions $c_{n}:(-1,1)\to \R^{+}$ are strictly monotonically 
increasing.
\end{cor}
\begin{proof} Since
    \[
    c_{n}(z) = \frac{\ds 1}{\ds n!}\vphi^{(n)}(z),
    \]
    it follows that $c_{n}'(z) = (n+1)c_{n+1}(z)>0$
\end{proof}

Furthermore, it is possible to show that for $z_{0}$ on $[0,1)$, 
these coefficients are also increasing functions of $n$, while for
$z_{0}$ on $[-1,-1/2)$ they will be decreasing.
\begin{thm}
    We have that
    \[
    c_{n}(z_{0})-c_{n-1}(z_{0})\left\{
    \begin{array}{ll}
    >0, & z_{0}\in[0,1)\\
    <0, & z_{0}\in[-1,-1/2]
    \end{array}
    \right.
    \]
for $n=2,3,\ldots$.
\end{thm}
\begin{proof}
    Let
    \[
    \theta(z) = \left[1-(z-z_{0})\right]\xi\left(\frac{\ds 1}{\ds 
    1-z}\right) = c_{0}+\dsum_{n=1}^{\infty}(c_{n}-c_{n-1})(z-z_{0})^n.
    \]
    and differentiate $\theta$ twice with respect to $z$ to obtain
    \[
    \theta''(z) = \frac{\ds 1+2z_{0}}{\ds (1-z)^3}\xi'\left(\frac{\ds 1}{\ds 
    1-z}\right)+\frac{\ds z_{0}}{(1-z)^4}\xi''\left(\frac{\ds 1}{\ds 
    1-z}\right).
    \]
    From this we see that $\theta''(z_{0})$ is positive for 
    $z_{0}\in[0,1)$
    and negative for $z_{0}\in[-1,-1/2]$. Further differentiation 
    with respect to $z$ gives rise to a sum of terms which are derivatives 
    of $\xi$ multiplied by constants which are either all positive or 
    all negative, depending on whether  $z_{0}\in[0,1)$ or 
    $z_{0}\in[-1,-1/2]$, respectively. Since in either case 
    $1/(1-z_{0})$ is greater than or equal to $1/2$, the derivatives 
    of $\xi$ are non--negative and at least one will be positive.
    
    The result now follows, as $c_{n}-c_{n-1} = \theta^{(n)}(z_{0})$ 
    for $n=2,3,\ldots$.
\end{proof}

Since the radius of convergence of the series for $\vphi$ around 
$z_{0}$ is 
$|1-z_{0}|$, it follows that, for negative $z_{0}$, $c_{n}(z_{0})\to0$ as
$n$ goes to infinity. We shall now prove that for all $z_{0}$ on 
$[0,1)$ these coefficients must go to infinity as $n$ goes to infinity.
    
\begin{thm} For $z_{0}\in[0,1)$ we have that
    \[
    \lim_{n\to \infty} c_{n}(z_{0}) = \infty.
    \]
\end{thm}
\begin{proof}
    From the expansion for $\xi$ given by~(\ref{xiseries}) it follows that
    \[
    \begin{array}{lll}
    \vphi(z) & = & \dsum_{n=0}^{\infty} s_{n}(1) \left(\frac{\ds z}{\ds 
    1-z}\right)^n\eqskip\\
    & = & s_{0}(1) + \dsum_{n=0}^{\infty} \left[\dsum_{k=0}^{n}
    \left(\begin{array}{c}n\\ 
    k\end{array}\right) s_{k+1}(1)\right] z^{n+1}.
    \end{array}
    \]
    This yields
    \[
    c_{n+1}(0) = \dsum_{k=0}^{n} \left(\begin{array}{c}n\\ 
    k\end{array}\right) s_{k+1}(1)\geq 
    s_{1}(1)+n\dsum_{k=1}^{n}\frac{\ds 1}{\ds k!} s_{k+1}(1),
    \]
    from which the result follows, since the $s_{k}$'s are positive and,
    by Corollary~\ref{monotz}, the coefficients $c_{n}$ are 
    increasing with $z$.
\end{proof}

\section{\label{proofmain}Proof of Theorem~\ref{mainthm}}

We shall now turn to the series for the 
function $\psi$ and the proof of Theorem~\ref{mainthm}. In order to 
do this, we begin by establishing the equivalence of the different 
formulations for the functions $\alpha_{n}(\tau)$.

\begin{lem}\label{lemequiv}
    For each positive integer $n$ and all real $\tau$ the functions 
    $\alpha_{n}(\tau)$ defined by~(\ref{dcoeff}) satisfy the 
    following relations:
    \begin{itemize}
    \item[(i)]
    \[
    \alpha_{n}(\tau) = \frac{\ds 1}{\ds \tau}
    \dsum_{\rho}\left[1-\left(\frac{\ds \rho}
    {\ds \rho-\tau}\right)^{n}\right],
    \]
    where $\rho$ runs over the nontrivial zeros of the zeta function 
    and the terms corresponding to $\rho$ and $1-\rho$ are  
    paired together.
    \item[(ii)] Define the power series coefficients $d_{n}(z_{0})$ by
    \[
    \psi(z) = \frac{\ds \vphi'(z)}{\ds \vphi(z)} =
    \dsum_{n=0}^{\infty} d_{n}(z_{0}) (z-z_{0})^n.
    \]
    Then
    \[
    \alpha_{n}(\tau) = \frac{\ds 1}{\ds 
    \tau^{n+1}}d_{n-1}(1-\frac{\ds 1}{\ds \tau})\;\; (\tau\neq0).
    \]
    \end{itemize}
\end{lem}
\begin{proof}

Starting from the product representation of $\xi$, that is,
\[
\xi(s) = \prod_{\rho}\left( 1 - \frac{\ds s}{\rho}\right),
\]
where, as usual, the product runs over the nontrivial zeros of the 
zeta function and each term is paired with that corresponding to 
$1-\rho$, we obtain that
\[
\vphi(s) = \prod_{\rho}\left(\frac{\ds 1-1/\rho-s}{\ds 1-s}\right).
\]
By taking the logarithmic derivative, 
we then have that
\[
\begin{array}{lll}
\psi(s) & = & -\dsum_{\rho}\left(\frac{\ds 1}{\ds \rho-1-\rho s}
\frac{\ds 1}{\ds 
1-s}\right)\eqskip\\
& = & -\dsum_{\rho}\frac{\ds 1}{\ds \rho(1-z_{0})-1}
\frac{\ds 1}{\ds 1-\frac{\ds \rho(s-z_{0})}{\ds \rho(1-z_{0})-1}}
\frac{\ds 1}{\ds 1-z_{0}}
\frac{\ds 1}{\ds 1-\frac{\ds s-z_{0}}{\ds 1-z_{0}}}\eqskip\\
& = & \dsum_{n=0}^{\infty}\tau^{n+1}\dsum_{\rho}\left[1-\left(\frac{\ds 
\rho}{\ds \rho-\tau}\right)^{n+1}\right](s-z_{0})^n,
\end{array}
\]
where $z_{0}$ is any real number on $[-1,1)$ and $\tau = 1/(1-z_{0})$.
We thus obtain that
\begin{equation}\label{dcoeffeq}
d_{n}(z_{0}) = \tau^{n+1}\dsum_{\rho}\left[1-\left(\frac{\ds 
\rho}{\ds \rho-\tau}\right)^{n+1}\right].
\end{equation}

On the other hand, we have that
\begin{equation}\label{alphaexpderiv}
\begin{array}{lll}
\alpha_{n+1}(\tau) &=& \frac{\ds 1}{\ds n!}\frac{\ds d^{n+1}}{\ds d 
s^{n+1}}\left[ s^{n}
\log\left(\xi(s)\right)\right]_{s=\tau}\eqskip\\ & = & -\frac{\ds 1}{\ds 
\tau}
\dsum_{\rho}\dsum_{k=0}^{n}\left(\begin{array}{c}n+1\\ 
    k\end{array}\right)\left(\frac{\ds \rho}{\ds 
    \tau}-1\right)^{k-n-1}\eqskip\\
    & = & \frac{\ds 1}{\ds \tau}\dsum_{\rho}\left[1-\left(\frac{\ds \rho}
    {\ds \rho-\tau}\right)^{n+1}\right].
\end{array}
\end{equation}
which proves $(i)$, and $(ii)$ now follows by comparing this with the 
expression obtained for $d_{n}(z_{0})$. 
\end{proof}

To prove Theorem~\ref{mainthm}, we begin by noticing that
equation~(\ref{dcoeffeq}) yields that the coefficients $d_{n}(z_{0})$ 
will be non--negative if
\[
\left| \frac{\ds 
\rho}{\ds \rho-\tau}\right|\leq 1.
\]
Hence, if all zeros have real part smaller than or equal 
to $\tau/2$, part $(ii)$ of Lemma~\ref{lemequiv} implies that
the coefficients $\alpha_{n}(\tau)$ will be non--negative.

Proceeding in the same fashion as in Li's paper, we shall now obtain
a recurrence relation between the 
coefficients $c_{n}$ and $d_{n}$ at a point $z_{0}$ on 
$[-1,1)$ -- note that since $\vphi$ is positive on this interval, $\psi$ 
is analytic in a neighbourhood of any point on the same interval. 
More precisely, from the definition of $\psi$ we have that
\[
\dsum_{n=1}^{\infty}n c_{n}(z-z_{0})^{n-1} =
\left(\dsum_{n=0}^{\infty}c_{n}(z-z_{0})^{n}\right)
\left(\dsum_{n=0}^{\infty}d_{n}(z-z_{0})^{n}\right),
\]
yielding the relations
\begin{equation}\label{recrel}
    d_{n-1} = n\frac{\ds c_{n}}{\ds c_{0}} - \frac{\ds 1}{\ds c_{0}}
    \dsum_{k=1}^{n-1} c_{k} d_{n-k-1}, \;\; n=1,2,\ldots.
\end{equation}
From Corollary~\ref{cpos} we know that the coefficients $c_{n}$ are 
all positive. Still following the steps of Li's proof, assume now that the
coefficients $d_{n}$ are all non--negative. Then, $d_{n-1}\leq
nc_{n}/c_{0}$, and
\[
\dsum_{n=1}^{\infty}\left| d_{n-1}(z-z_{0})^{n-1}\right|\leq
\frac{\ds 1}{\ds c_{0}}\dsum_{n=1}^{\infty}n 
c_{n}\left|z-z_{0}\right|^{n-1} = \frac{\ds 1}{\ds c_{0}}
\vphi'(z_{0}+\left|z-z_{0}\right|).
\]
It follows that the series
\[
\dsum_{n=0}^{\infty}\left| d_{n}(z-z_{0})^{n}\right|
\]
converges for all $\left| z-z_{0}\right|<1-z_{0}$. This finishes the 
proof of Theorem~\ref{mainthm}.

\section{Some properties of the coefficients $\alpha_{n}$\label{prop}}

We shall now study the behaviour of the functions $\alpha_{n}$. 
Although Theorem~\ref{mainthm} only considers the case of $\tau$ 
larger than $1/2$, we will consider the whole of $\R^{+}$ unless it is 
explicitly stated otherwise. We shall stick to the notation used 
in the previous section.

\subsection{A relaxed sufficient condition}

We begin by giving a result in the spirit of Corollary 1 
in~\cite{bola}, namely, we show that the sufficient condition of 
Theorem~\ref{mainthm} can be somewhat relaxed.
\begin{prop} Let\label{mainprop}
    \[
    g(z) =\dsum_{n=0}^{\infty} \gamma_{n}(z-z_{0})^n
    \]
    be any function such that the coefficients $\gamma_{n}=\gamma_{n}(\tau)$ are
    all nonnegative and such that the Taylor series for $g$ around 
    $z_{0}$ has radius of convergence at least $1-z_{0}$. Then,
    given $\tau$ on $[1/2,\infty)$, the half--plane
    $\Re(s)>\tau/2$ is a zero--free region for the Riemann zeta
    function if and only if $\alpha_{n}(\tau)+\gamma_{n-1}(\tau)$ is
    nonnegative for all positive integer $n$.
\end{prop}
\begin{proof}
    That the RH implies $\alpha_{n}(\tau)+\gamma_{n-1}(\tau)\geq0$
is a trivial consequence of Theorem~\ref{mainthm} and the hypothesis that 
$\gamma_{n}\geq0$ for all $n$.

To prove the implication in the other direction we consider the 
auxiliary function
\[
h(z) = e^{\int_{0}^{z} g(t)dt}\vphi(z).
\]
We have that the Taylor series for $h$ around a point $z_{0}$ on $[-1,1)$
will also have radius of convergence $1-z_{0}$, and the corresponding 
Taylor coefficients are also nonnegative. It folows that
\[
\begin{array}{lll}
H(z) & := & \frac{\ds h'(z)}{\ds h(z)}\eqskip\\
 & = & \psi(z) + g(z)\eqskip\\
 & = & \dsum_{n=0}^{\infty} 
 \left[d_{n}(z_{0})+\gamma_{n}(z_{0})\right] (z-z_{0})^{n}.
 \end{array}
 \]
 Since the radius of convergence of the series for $g$ is $1-z_{0}$, 
 it follows that the radius of convergence for the Taylor series of 
 $H$ around $z_{0}$ will still be that of the series for $\psi$. 
 Proceeding now as in the second part of the proof of 
 Theorem~\ref{mainthm} with $d_{n}(z_{0})+\gamma_{n}(z_{0})$ in place of
 $d_{n}(z_{0})$ we obtain the desired result.
\end{proof}

A straightforward consequence of this result is that for a value of $\tau$ larger 
than $1/2$ for which the half--plane $\Re(s)>\tau/2$ is not a zero--free region
of the zeta function there cannot exist only a finite number of the 
$\alpha_{n}$'s which are negative.

\begin{cor}\label{infneg}
    If the half--plane $\Re(s)>\tau_{0}/2$ is not a zero--free region
of the Riemann zeta function for some $\tau_{0}$ on $[1/2,2)$, then 
there exists a strictly increasing infinite sequence $n'$ such that 
$\alpha_{n'}(\tau_{0})$ is negative.
\end{cor}

\subsection{A system of equations satisfied by the functions $\alpha_{n}$}

In this section we obtain a set of differential equations which are satisfied by 
the functions $\alpha_{n}(\tau)$. These will be used to establish some of 
the properties of these functions in the remaining of the paper.

\begin{thm}\label{tsysdiff}
    The functions $\alpha_{n}:\R^{+}\to\R$ satisfy the following 
    (infinite) system of differential equations
    \begin{equation}\label{sysdiff}
    \frac{\ds \tau}{\ds n} \alpha_{n}'(\tau)+
    \frac{\ds n+1}{\ds n}\alpha_{n}(\tau)=
    \alpha_{n+1}(\tau), \;\; n= 1,2,\ldots.
    \end{equation}
\end{thm}
\begin{proof}
    By~(\ref{alphaexpderiv}) we have that
    \[
    \begin{array}{lll}
    \alpha_{n+1}(\tau) & = & \frac{\ds 1}{\ds \tau}\dsum_{\rho}\left[1-
    \left(\frac{\ds \rho}
    {\ds \rho-\tau}\right)^{n+1}\right]\eqskip\\
    & = & 
    \frac{\ds 1}{\ds \tau}\dsum_{\rho}\left[1-
    \left(\frac{\ds \rho}
    {\ds \rho-\tau}\right)^{n}\right]+
    \frac{\ds 1}{\ds \tau}\dsum_{\rho}\left(\frac{\ds \rho}{\ds 
    \rho-\tau}\right)^{n}\left[1-
    \frac{\ds \rho}
    {\ds \rho-\tau}\right]\eqskip\\
    & = & \alpha_{n}(\tau)-\dsum_{\rho}\frac{\ds \rho^{n}}{\ds 
    (\rho-\tau)^{n+1}}.
    \end{array}
    \]
    On the other hand, differentiating the above expression for 
    $\alpha_{n}$ with respect to $\tau$ gives
    \[
    \alpha_{n}'(\tau) = -\frac{\ds 1}{\ds \tau}\alpha_{n}(\tau)-
    \frac{\ds n}{\ds \tau}\dsum_{\rho}\frac{\ds 
    \rho^{n}}{(\rho-\tau)^{n+1}}.
    \]
    Combining these two expressions yields the desired result.
\end{proof}
    
System~(\ref{sysdiff}) has a singularity at 
$\tau=0$ which means that the conditions at this point given by 
Theorem~\ref{oscilthm} should be seen as {\it compatibility} 
conditions necessary to ensure that the solutions have a bounded 
derivative. Thus, in order to determine the solution of this system 
uniquely, one should add a condition at another point which may 
actually be the behaviour at infinity. In any case, since we are 
dealing with an infinite system, a more careful analysis would be 
required to ensure uniqueness under these circumstances. Since for our
present purposes the fact that the coefficients $\alpha_{n}$ 
satisfy~(\ref{sysdiff}) is sufficient, we will not pursue the matter 
here.

We also remark that although this system is nonautonomous, it can be 
transformed into an autonomous system via a change of variables. This 
is the case if we take, for instance, $\beta_{n}(\tau) = 
e^{-\tau}\alpha_{n}(e^{-\tau})$, yielding
\[
\beta_{n}'(\tau) - n\beta_{n}(\tau) = -n\beta_{n+1}(\tau).
\]
The interval of interest is now $[-\log 2,+\infty)$, and we have that
\[
\beta_{n}(-\log 2)= 2\alpha_{n}(2)  (\geq0) \mbox{ and }
\lim_{\tau\rightarrow +\infty}\beta_{n}(\tau) = 0^{-}.
\]

Another approach is to view system~(\ref{sysdiff}) as a discrete 
dynamical system defined on the set of 
$\mathcal{C}^{\infty}(0,\infty)$ functions, with initial 
condition $\alpha_{1}(\tau) = \xi'(\tau)/\xi(\tau)$. In this setting the 
problem is well--posed, but again we are dealing with a 
non--autonomous sytem.

\subsection{Oscillatory behaviour of the functions $\alpha_{n}$}

We shall now study the oscillatory properties of the functions 
$\alpha_{n}$. We begin by proving the following simple lemma.
\begin{lem}\label{intlem}
    For $0\leq \tau_{1}<\tau_{2}$ we have
    \[
    \int_{\tau_{1}}^{\tau_{2}}t^{n}\alpha_{n+1}(t)dt =
    \frac{\ds \tau_{2}^{n+1}\alpha_{n}(\tau_{2})-\tau_{1}^{n+1}\alpha_{n}
    (\tau_{1})}{\ds n}, \;\; n=1,\ldots.
    \]
\end{lem}
\begin{proof}
    Multiplying equation~(\ref{sysdiff}) by 
    $\tau^{n+1}$ we get
    \[
    \frac{\ds d}{\ds d\tau}\left[\tau^{n+1}\alpha_{n}(\tau)\right] =
    n\tau^{n}\alpha_{n+1}(\tau).
    \]
    Integrating now between $\tau_{1}$ and $\tau_{2}$ gives the desired 
    result.
\end{proof}

We are now ready to prove Theorem~\ref{oscilthm}.

Proceeding as in the proof of Theorem~\ref{tsysdiff}  and letting 
$\tau=0$ we get that
\[
\alpha_{n+1}(0) = \alpha_{n}(0)-\dsum_{\rho}\frac{\ds 1}{\ds 
   \rho} = \alpha_{n}(0)+\frac{\ds \xi'(0)}{\ds \xi(0)}.
\]
On the other hand, from the definition of the functions $\alpha_{n}$ 
we have that $\alpha_{1}(0) = \xi'(0)/\xi(0)$ and so $\alpha_{n}(0)=
n\xi'(0)/\xi(0)$. This proves part (i) of 
the theorem.

To prove (ii) we argue by induction. Since $\alpha_{1}(\tau) = 
\xi'(\tau)/\xi(\tau)$, we have that $\alpha_{1}(1/2)=0$. Letting 
$\tau_{1}=0$ and $\tau_{2}=1/2$ in Lemma~\ref{intlem} we obtain that
\[
\dint_{0}^{1/2} t \alpha_{2}(t)dt = 0
\]
and thus $\alpha_{2}$ must have a zero on $(0,1/2)$.
Assume now that $\alpha_{n}(a_{n})=0$ for some $a_{n}$ less than 
$1/2$. Then, letting 
$\tau_{1}=0$ and $\tau_{2}=a_{n}$ in Lemma~\ref{intlem} we obtain that
\[
\dint_{0}^{a_{n}} t^{n} \alpha_{n+1}(t)dt = 0,
\]
and the result follows.

Point (iii) in the theorem follows, in the case of distinct zeros
$z_{1}<z_{2}$ of $\alpha_{n}$, in a similar fashion using
now Lemma~\ref{intlem} with $\tau_{1}=z_{1}$ and $\tau_{2}=z_{2}$. 
For the case of a zero $z$ of multiplicity $m$ of $\alpha_{n}$, 
it follows by repeated differentiation of~(\ref{sysdiff}) that 
$\alpha_{n+1}$ vanishes to order $m-1$ at $z$.

In order to prove point (iv), denote by $\nu_{n}$ the number of zeros 
of $\alpha_{n}$, counting multiplicities. From 
points (ii) and (iii) it follows that $\nu_{n}\leq \nu_{n+1}$. We argue by 
contradiction. If the statement is false, there exists a smallest number 
$A$ such that $\nu(n)\leq A$ for all $n$. Let $p$ be the smallest number 
for which $\nu_{p}=A$, and denote the largest 
zero of $\alpha_{p}$ by $z_{p}$ and its multiplicity by $m$.
Then $\alpha_{p+m}$ has all its $A$ positive zeros strictly to the left 
of $z_{p}$. Denote its largest zero by $z_{p+m}$. Since for $n$ larger 
than $p+m$ all the zeros of the functions $\alpha_{n}$ must be to the 
left of $z_{p+m}$, it follows that there is only a finite number of 
the functions $\alpha_{n}$ which are negative on the interval 
$(z_{m+p},z_{p})$. Then  Proposition~\ref{mainprop} gives
that there are no zeros of $\xi$ in the half--plane $\Re(s)>z_{p+m}$.

On the other hand, we have that the zero at $z_{p}$ of at least one of the
functions $\alpha_{p},\ldots, \alpha_{p+m-1}$ is of odd multiplicity, 
which means that there is at least one function which takes on 
negative values on the interval $(z_{p+m},z_{p})$. By 
Theorem~\ref{mainthm}, this implies that there must exist zeros of the 
function $\xi$ with real parts between $z_{p+m}$ and $z_{p}$, giving a 
contradiction. This concludes the proof of Theorem~\ref{oscilthm}.

\section{A function defined on $\R^{2}$\label{secf}}

In this section we extend the definition of the coefficients 
$\alpha_{n}$ to the case where the parameter $n$ is allowed to vary
continuously. This is done by means of one of the expressions for 
$\alpha_{n}$ used in the previous section, namely, identity~(\ref{alphaexpderiv}).
More precisely, define the function
$F:\R^{2}\to\R$ by
\[
F(x,\tau) = \frac{\ds 1}{\ds \tau}\sum_{\rho}\left[1-\left(\frac{\ds\rho}
{\ds\rho-\tau}\right)^{x}\right].
\]
From this and the definition of the coefficients 
$\alpha_{n}$ it follows that $\alpha_{n}(\tau)=F(n,\tau)$ for 
$n=1,\ldots$. Many of the properties derived in the previous section 
also apply to the function $F$. For instance, it is possible to prove 
in the same way that it will also satify a 
corresponding system of differential equations, namely,
\[
\frac{\ds \tau}{\ds x}\frac{\partial F(x,\tau)}{\ds \partial\tau}+
\frac{\ds x+1}{\ds x}F(x,\tau)=F(x+1,\tau).
\]
It is also clear that $F$ vanishes identically in $\tau$ at $x=0$ and taking 
derivatives with respect to $x$ and letting $x=0$ yields
\[
\frac{\ds \partial F(0,\tau)}{\ds \partial x} = \frac{\ds 1}{\ds \tau}
\log\left[\xi(\tau)\right],
\]
giving that this derivative is negative for $\tau$ smaller than one 
and positive for $\tau$ larger than one.

\subsection{Alternative representations of the function $F$}

There are several ways in which one can transform $F$ in 
order to obtain other representations of this function. Here we give 
two which depend on the Taylor coefficients of the logarithmic derivative of 
$\xi$.

\begin{thm}\label{altrepf} We have that
    \[
    F(x,\tau) = \dsum_{k=0}^{\infty} \frac{\ds \ell_{k}}{\ds (k+1)!}
    \frac{\ds \Gamma(x+k+1)}{\ds \Gamma(x)}\tau^{k},\;\; |\tau|<R, \;
    x\in\CM,
    \]
    where $\ell_{k}$ are the Taylor coefficients of the function 
    $\ell(s) = \xi'(s)/\xi(s)$ around zero and $R$ is the absolute value of
    the first zero of the zeta function on the critical line.
\end{thm}
\begin{proof}
    Write
    \[
    G(x,\tau) = \tau F(x,\tau) = \dsum_{\rho}\left[1-\left(\frac{\ds\rho}
    {\ds\rho-\tau}\right)^{x}\right].
    \]
    Then\[
    \begin{array}{lll}
    \frac{\ds \partial^{k} G(x,0)}{\ds \partial\tau^{k}} & = &
    -x(x+1)\ldots(x+k-1)\dsum_{\rho}\frac{\ds 1}{\ds \rho^{k}} \eqskip\\
     & = &
    x(x+1)\ldots(x+k-1)\frac{\ds 1}{\ds (k-1)!}\ell^{(k-1)}(0),
    \end{array}
    \]
    and
    \[
    G(x,\tau) = \dsum_{k=1}^{\infty}\frac{\ds 1}{\ds k!} \frac{\ds 
    \Gamma(x+k)}{\ds \Gamma(x)}\ell_{k-1}\tau^{k},
    \]
    from which the result follows for $F$ upon multiplication by 
    $\tau$ and a change of the summation variable.
    
    The value of the radius of convergence follows from the fact that the
    only singularities of $F$ (in $\tau$) are at the zeros of $\xi$.
\end{proof}

From this it follows that at negative integers $F$
is a polynomial in $\tau$ -- this could also have been obtained by 
direct inspection of the definition of $F$. More precisely,
we have the following
\begin{cor}
    Let $p$ be a positive integer. Then
    \[
    F(-p,\tau) = \dsum_{k=0}^{p-1} \frac{\ds \ell_{k}}{\ds 
    (k+1)!}\prod_{j=0}^{k}(j-p)\tau^{k}.
    \]
\end{cor}

We also obtain a generalization of an expression given in~\cite{coff2} 
for Li's coefficients, which is the sum of a positive and a negative 
part. To this end, we need to define the coefficients $\eta_{k}$ 
which appear in the  Taylor series development of the logarithmic 
derivative of $\zeta$ around one. More precisely,
\[
\frac{\ds \zeta'(s)}{\ds \zeta(s)} = -\frac{\ds 1}{\ds s-1}-
\dsum_{k=0}^{\infty}\eta_{k}(s-1)^{k}.
\]
These coefficients may also be written in terms of the von Mangoldt 
function~--~see~\cite{coff2}, for instance.

\begin{cor} Define
    \[
    P(x,\tau) = \dsum_{k=1}^{\infty}\frac{\ds 1}{\ds (k+1)!}
    \frac{\ds \Gamma(x+k+1)}{\ds 
    \Gamma(x)}\left[\left(1-2^{(-k-1)}\right)\zeta(k+1)-1\right]\tau^{k}
    \]
    and
    \[
    N(x,\tau) = \left[\log\left(2\sqrt{\pi}\right)-1-\frac{\ds \gamma}{\ds 2}
    \right]x+\dsum_{k=1}^{\infty}\frac{\ds (-1)^{k}}{\ds (k+1)!}
    \frac{\ds \Gamma(x+k+1)}{\ds 
    \Gamma(x)}\eta_{k}\tau^{k}.
    \]
    Then, for positive $x$ and $\tau$, the functions $P$ and $N$ are 
    positive and negative, respectively, and    
    \[
    F(x,\tau) = P(x,\tau) + N(x,\tau).
    \]
\end{cor}
\begin{rem} This result generalizes one of the expressions for Li's 
coefficients given in~\cite{coff2}, while explicitly separating
the positive and negative parts of the expansion.
\end{rem}
\begin{proof}
    We have that
    \[
    \log[\xi(s)] = \log\left(\frac{\ds s}{\ds 2}\right)-\frac{\ds s}{\ds 2}
    \log(\pi)+
    \log\left[\Gamma\left(\frac{\ds s}{\ds 2}\right)\right]-\dsum_{k=1}^{\infty}
    \frac{\ds\eta_{k-1}}{\ds k}(s-1)^{k},
    \]
    from which it follows upon exchanging $s$ with $1-s$ and 
    differentiating with respect to $s$ that
    \[
    \ell(s) = -\frac{\ds 1}{\ds 1-s}+\log\left(\sqrt{\pi}\right)-
    \frac{\ds 1}{\ds 2}\psi\left(\frac{\ds 1-s}{\ds 2}\right)+
    \dsum_{k=0}^{\infty}(-1)^k\eta_{k}s^{k},
    \]
    where $\psi$ is the digamma function defined by $\psi(s) 
    =\left[\log\left(\Gamma(s)\right)\right]'$. Since
    \[
    -\frac{\ds 1}{\ds 2}\psi\left(\frac{\ds 1-s}{\ds 2}\right) =
    \frac{\ds \gamma}{\ds 2}+\log(2)+\dsum_{k=1}^{\infty}\left(
    1-2^{-k-1}\right)\zeta(k+1)s^{k},
    \]
    we are finally led to
    \[
    \ell(s) = \log\left(2\sqrt{\pi}\right)-1-\frac{\ds \gamma}{\ds 2}+
    \dsum_{k=1}^{\infty}\left[-1+(1-2^{-k-1})\zeta(k+1)+(-1)^{k}\eta_{k}
    \right]s^{k}.
    \]
    This last expression gives the coefficients $\ell_{k}$ used in 
    Theorem~\ref{altrepf}, from which the expression for $F$ in terms 
    of $P$ and $N$ follows. To prove that $P$ and $N$ are positive end 
    negative, respectively, it is sufficient to note that 
    $(1-2^{(-k-1)})\zeta(k+1)-1$ is positive, and that the 
    coefficients $\eta_{k}$ are positive for odd $k$ and negative for 
    even $k$. The first of these statements follows directly from the 
    fact that $\zeta(k+1)$ is larger than $1+2^{-(k+1)}+3^{-(k+1)}$, while 
    the second can be found in~\cite{coff2}.    
\end{proof}

We shall now show the existence of a region close to the $x$ axis 
where the function $F$ remains negative. In particular, this will 
imply that the sequence $a_{n}$ referred to in Theorem~\ref{oscilthm},
if it goes to zero, cannot do so faster than $c/n$, for some positive 
constant $c$.

\begin{thm} 
    There exist positive constants $\tau_{0}$ and $x_{0}$ such that 
    $F$ takes on negative values on the set
    \[
    X=\left\{(x,\tau)\in\R^{2}: x_{0}<x<\tau_{0}/\tau \right\}.
    \]
\end{thm}
\begin{proof}
    From Theorem~\ref{altrepf}, we have that
    \[
    r F(\frac{\ds 1}{\ds r},r\tau) = 
    \dsum_{k=0}^{\infty}\frac{\ds \ell_{k}}{\ds (k+1)!}
    \left( 1+r\right)\ldots\left( 1+kr\right)\tau^{k}
    \]
    and thus
    \begin{equation}\label{limr}
    \lim_{r\to 0}\left[ rF(\frac{\ds 1}{\ds r},r\tau)\right] = 
    \dsum_{k=0}^{\infty} \frac{\ds \ell_{k}}{\ds 
    (k+1)!}\tau^k.
    \end{equation}       
    Since $\ell_{0}$ is negative, there exists a positive number 
    $\tau_{0}$ such that the function appearing on the right--hand 
    side of~(\ref{limr}) is negative at $\tau_{0}$. Hence, for 
    positive $r$ sufficiently small, $F(1/r, rt)$ must be negative 
    for $0<t<\tau_{0}$. To obtain the result, write $x=1/r$ and $\tau=rt$.
\end{proof}

The above result does not show if the region close to the $x-$axis 
thins out as $x$ goes to infinity, or if there is a strip where $F$ is 
negative. In fact, we conjecture that the latter situation does not 
occur, and thus that the sequence 
$a_{n}$ in Theorem~\ref{oscilthm} does indeed go to zero.

If we use the binomial series development for the power term in the 
definition of $F$ we get a slightly different expression.
\begin{thm}
    We have that
    \[
    F(x,\tau) = \dsum_{k=0}^{\infty}\frac{\ds \Gamma(x+1)}{\ds 
    (k+1)!\Gamma(x-k)}\ell_{k}(\tau)\tau^{k},
    \]
    where
    \[
    \ell_{k}(\tau) = \frac{\ds 1}{\ds k!}\frac{\ds d^{k}}{\ds 
    d\tau^{k}}\left[\ell(\tau)\right].
    \]
\end{thm}
\begin{proof}
    \[
    \begin{array}{lll}
    F(x,\tau) & = & \frac{\ds 1}{\ds \tau}\dsum_{\rho}\left[1-\left(\frac{\ds\rho}
    {\ds\rho-\tau}\right)^{x}\right]\eqskip\\
    & = & \frac{\ds 1}{\ds 
    \tau}\dsum_{\rho}\left[1-\left(1+\frac{\ds\tau}
    {\ds\rho-\tau}\right)^{x}\right]\eqskip\\
    & = & \frac{\ds 1}{\ds 
    \tau}\dsum_{\rho}\left[1-\dsum_{k=0}^{\infty}\frac{\ds 
    \Gamma(x+1)}{\ds k!\Gamma(x+1-k)}\left(\frac{\ds \tau}{\ds 
    \rho-\tau}\right)^{k}\right]\eqskip\\
    & = & -\dsum_{k=1}^{\infty}\frac{\ds 
    \Gamma(x+1)}{\ds k!\Gamma(x+1-k)} \tau^{k-1}\dsum_{\rho}\frac{\ds 
    1}{\ds (\rho-\tau)^k}\eqskip\\
    & = & \dsum_{k=1}^{\infty}\frac{\ds 
    \Gamma(x+1)}{\ds k!\Gamma(x+1-k)} \frac{\ds \tau^{k-1}}{\ds 
    (k-1)!} \frac{\ds d^{k-1}}{\ds 
    d\tau^{k-1}}\left[\ell(\tau)\right],
    \end{array}
    \]
    which gives the desired result.
\end{proof}

\section{Concluding remarks}

As was pointed out in the Introduction, the criterion
given here is not restricted to
the $\xi$ function, and most of the results presented may be extended to a 
fairly general class of complex functions. In particular, we mention 
that the set of equations~(\ref{sysdiff}) satisfied by the 
coefficients $\alpha_{n}$ (and the function $F$) is independent of 
the function under study, and that the sequence 
$\alpha_{n}$ will be determined by the function $\alpha_{1}$, for 
instance. The alternative representations for the function $F$ 
given in Section~\ref{secf} also apply to this more general case, if 
we replace the Taylor coefficients of the function $\ell$ by those of 
the logarithmic derivative of a function $f$. In 
relation to this, we stress again that for a function $f$ of the
form~(\ref{prodf}) it would also have been possible to use the more general
approach of Bombieri and Lagarias~\cite{bola} throughout the paper, in which
case it would not have been necessary to require the non--negativity
condition on the Taylor coefficients of $f$.

Regarding the specific case of the function $\xi$, we note that for
$z_{0}$ in $(-1,0)$ the disk $\left| z-z_{0}\right|<1-z_{0}$
contains the closure of the unit disk centred at
zero, with the exception of the point $z=1$. Since $\vphi$ 
vanishes at the points on the unit circle which correspond to the 
nontrivial zeros of $\zeta$ on the critical line, $\psi$ cannot be 
analytic on the larger circle which implies that the coefficients 
$d_{n}(z_{0})$ must then take on negative values when $z_{0}$ is negative,
and hence when $\tau$ is less than one we have the same type of 
behaviour for the coefficients $\alpha_{n}$. 

On the other hand, and as was also pointed out in the Introduction, for 
$\tau$ greater than or equal to two, the half--plane $\Re(z)<\tau/2$
contains the whole of the critical strip, and hence the $\alpha_{n}$'s
will be non--negative in this case.

Besides giving a criterion for the existence/nonexistence of zeros 
in half--planes, we hope that the results presented here will also 
help to provide some insight into the behaviour of the coefficients
$\lambda_{n}$. In particular, it might be interesting to analyse 
numerically the situation where $\tau$ is close to but less than one, 
and see how the negative part of the sequence $\alpha_{n}$ behaves as $\tau$ 
approaches one. Note also that when $\tau$ is zero, all functions 
$\alpha_{n}$ (with positive $n$) are negative.

\subsection*{Acknowledgments} I would like to thank Mark Coffey for having 
made available a preliminary version of~\cite{coff2}.

\end{document}